\magnification=1100
\def\to{\ \longrightarrow\ }

\def\nl{\hfill\break}

\def\hexnumber#1{\ifcase#1 0\or 1\or 2\or 3\or 4\or 5\or 6\or 7\or 8\or
 9\or A\or B\or C\or D\or E\or F\fi}
%
%
\font\twelvemsa=msam10 scaled 1200   
\font\tenmsa=msam10                  
\font\ninemsa=msam9            \font\sevenmsa=msam7
\font\sixmsa=msam6             \font\fivemsa=msam5
%
%
\newfam\msafam                 \textfont\msafam=\tenmsa
\scriptfont\msafam=\sevenmsa   \scriptscriptfont\msafam=\fivemsa
\edef\hexa{\hexnumber\msafam}        
\def\msa{\fam\msafam\tenmsa}         
%
%
\font\twelvemsb=msbm10 scaled 1200   
\font\tenmsb=msbm10                  
\font\ninemsb=msbm9            \font\sevenmsb=msbm7
\font\sixmsb=msbm6             \font\fivemsb=msbm5
%
\newfam\msbfam                 \textfont\msbfam=\tenmsb       
\scriptfont\msbfam=\sevenmsb   \scriptscriptfont\msbfam=\fivemsb
\edef\hexb{\hexnumber\msbfam}        
\def\msb{\fam\msbfam\tenmsb}         
%
%
\font\twelveeufm=eufm10 scaled 1200  
\font\teneufm=eufm10                 
\font\nineeufm=eufm9           \font\seveneufm=eufm7
\font\sixeufm=eufm6            \font\fiveeufm=eufm5
%
\newfam\eufmfam                \textfont\eufmfam=\teneufm
\scriptfont\eufmfam=\seveneufm \scriptscriptfont\eufmfam=\fiveeufm
\edef\hexf{\hexnumber\eufmfam}      
\def\frak{\fam\eufmfam\teneufm}     
%
%
%
\font\twelverm=cmr10 scaled 1200    
\font\ninerm=cmr9                   
\font\sixrm=cmr6   
%
\font\twelvei=cmmi10 scaled 1200    
\font\ninei=cmmi9                   
\font\sixi=cmmi6  
%
\font\twelvesy=cmsy10 scaled 1200   
\font\ninesy=cmsy9                  
\font\sixsy=cmsy6  
%
\font\twelvebf=cmbx10 scaled 1200   
\font\ninebf=cmbx9                  
\font\sixbf=cmbx6  
%
%
\font\twelveit=cmti10 scaled 1200   
\font\nineit=cmti9                  
%
\font\twelvesl=cmsl10 scaled 1200   
\font\ninesl=cmsl9                  
%
\font\twelvett=cmtt10 scaled 1200   
\font\ninett=cmtt9                  
%
%
%
%
\def\small{%
%
%
\textfont0=\ninerm \scriptfont0=\sixrm \scriptscriptfont0=\fiverm
\def\rm{\fam0\ninerm}        
%
%
\textfont1=\ninei \scriptfont1=\sixi \scriptscriptfont1=\fivei
%
%
\textfont2=\ninesy \scriptfont2=\sixsy \scriptscriptfont2=\fivesy
%
%
\textfont3=\tenex \scriptfont3=\tenex \scriptscriptfont3=\tenex
%
%
\textfont\bffam=\ninebf \scriptfont\bffam=\sixbf
\scriptscriptfont\bffam=\fivebf \def\bf{\fam\bffam\ninebf}%
%
%
\textfont\itfam=\nineit \def\it{\fam\itfam\nineit}%
\textfont\slfam=\ninesl \def\sl{\fam\slfam\ninesl}%
\textfont\ttfam=\ninett \def\tt{\fam\ttfam\ninett}%
%
%
%
\textfont\msafam=\ninemsa \scriptfont\msafam=\sixmsa
\scriptscriptfont\msafam=\fivemsa \def\msa{\fam\msafam\ninemsa}%
%
%
\textfont\msbfam=\ninemsb \scriptfont\msbfam=\sixmsb
\scriptscriptfont\msbfam=\fivemsb \def\msb{\fam\msbfam\ninemsb}%
%
%
\textfont\eufmfam=\nineeufm  \scriptfont\eufmfam=\sixeufm
\scriptscriptfont\eufmfam=\fiveeufm \def\frak{\fam\eufmfam\nineeufm}%
%
%
%
\normalbaselineskip=11pt
\setbox\strutbox=\hbox{\vrule height8pt depth3pt width0pt}%
%
%
\normalbaselines\rm}    
%
%
%
%
\def\large{%
\textfont0=\twelverm \scriptfont0=\ninerm \scriptscriptfont0=\sevenrm
\def\rm{\fam0\twelverm}%
\textfont1=\twelvei \scriptfont1=\ninei \scriptscriptfont1=\seveni
\textfont2=\twelvesy \scriptfont2=\ninesy \scriptscriptfont2=\sevensy
\textfont3=\tenex \scriptfont3=\tenex \scriptscriptfont3=\tenex
\textfont\bffam=\twelvebf \scriptfont\bffam=\ninebf
\scriptscriptfont\bffam=\sevenbf \def\bf{\fam\bffam\twelvebf}%
\textfont\itfam=\twelveit \def\it{\fam\itfam\twelveit}%
\textfont\slfam=\twelvesl \def\sl{\fam\slfam\twelvesl}%
\textfont\ttfam=\twelvett \def\tt{\fam\ttfam\twelvett}%
\textfont\msafam=\twelvemsa \scriptfont\msafam=\ninemsa
\scriptscriptfont\msafam=\sevenmsa \def\msa{\fam\msafam\twelvemsa}         
\textfont\msbfam=\twelvemsb \scriptfont\msbfam=\ninemsb
\scriptscriptfont\msbfam=\sevenmsb \def\msb{\fam\msbfam\twelvemsb}         
\textfont\eufmfam=\twelveeufm  \scriptfont\eufmfam=\nineeufm
\scriptscriptfont\eufmfam=\seveneufm \def\frak{\fam\eufmfam\teneufm}
\normalbaselineskip=15pt
\setbox\strutbox=\hbox{\vrule height11pt depth4pt width0pt}%
\normalbaselines\rm}%
%
\def\Bbb{\msb}

%

%
\mathchardef\plussquare="0\hexa01
\mathchardef\nge="3\hexb0B
\mathchardef\maltesecross="0\hexa7A
\mathchardef\del="0\hexf01
%

%

 \input epsf
\overfullrule=0pt
\mathsurround=2pt

\font\Bbb=msbm10

\font\secfont=cmbx10
\font\ab=cmbx8

\font\nam=cmr8
\font\aff=cmti8

\font\teneufm=eufm10

\mathchardef\square="0\hexa03
\def\qed{\hfill$\square$\par\rm}
\def\np{\vfill\eject}
\def\boxing#1{\ \lower 3.5pt\vbox{\vskip 3.5pt\hrule \hbox{\strut\vrule
\ #1 \vrule} \hrule} }

\def\down#1{\ \lower 3.5pt\vbox{\vskip 3.5pt \hbox{\strut \ #1 \vrule} \hrule}
}
\def\negdown#1{\ \lower 3.5pt\vbox{\vskip 3.5pt \hbox{\strut  \vrule \ #1
}\hrule} }

\hsize=6.3 truein
\vsize=9 truein

\baselineskip=13 pt
\parskip=\baselineskip
 1

\parindent=0pt

\def\Z{\hbox{\Bbb Z}}
\def\R{\hbox{\Bbb R}}

\def\H{\hbox{\Bbb H}}

\def\I{\hbox{\Bbb I}}

\def\s{\sigma}
\def\t{\tau}






\newif \iftitlepage \titlepagetrue

\def\diagram{\global\advance\diagramnumber by 1
$$\epsfbox{turfig.\number\diagramnumber}$$}
\def\ddiagram{\global\advance\diagramnumber by 1
\epsfbox{turfig.\number\diagramnumber}}

\newcount\diagramnumber
\diagramnumber=0

\newcount\secnum \secnum=0
\newcount\subsecnum
\newcount\defnum
\def\section#1{
               \vskip 10 pt
               \advance\secnum by 1 \subsecnum=0
               \leftline{\secfont \the\secnum \quad#1}
               }

\def\subsection#1{
               \vskip 10 pt
               \advance\subsecnum by 1
               \defnum=1
               \leftline{\secfont \the\secnum.\the\subsecnum\ \quad #1}
               }

\def\definition{
               \advance\defnum by 1
               \bf Definition
\the\secnum .\the\defnum \rm \
               }

\def\lemma{
               \advance\defnum by 1
               \par\bf Lemma  \the\secnum
.\the\defnum \sl \ \par \rm
               }

\def\theorem{
               \advance\defnum by 1
               \par\bf Theorem  \the\secnum
.\the\defnum \rm \ 
              }

\def\cite#1{\secfont [#1]\rm}

\vglue 20 pt

\centerline{\secfont Biquandles of Small Size and some Invariants of Virtual and Welded Knots}

\medskip

\centerline{\nam ANDREW BARTHOLOMEW, ROGER FENN}
\centerline{\aff School of Mathematical Sciences, University of Sussex}
\centerline{\aff Falmer, Brighton, BN1 9RH, England}

\bigskip

\baselineskip=10 pt
\parskip=0 pt
\bigskip
\centerline{\nam ABSTRACT}
\leftskip=0.25 in
\rightskip=0.25in

{\ab In this paper we give the results of a computer search for biracks of small size and we give various interpretations of these findings. The list includes biquandles, racks and quandles together with new invariants of welded knots and examples of welded knots which are shown to be non-trivial by the new invariants. These can be used to answer various questions concerning virtual and welded knots. As an application we reprove the result that the Burau map from braids to matrices is non injective and give an example of a non-trivial virtual (welded) knot which cannot be distinguished from the unknot by any linear biquandles.}

\medskip
\leftskip=0 in
\rightskip=0 in
\baselineskip=13 pt
\parskip=\baselineskip

\centerline{{\bf Key words:}  birack, biquandle, rack, quandle, virtual and welded knot}

\section{Introduction}
This paper is based on a computer search for biracks, biquandles, racks and quandles, \cite{FJK, K, FR, J, M}, of small size. The search is complete for size up to 4 and agrees with other searches, \cite{N}, where they overlap. Size $n=5$ seems to be difficult but it may succumb with a different programming technique. However for $n=5, 6$ we produce a complete list of racks and this allows us to produce a complete list of so called rack-related biracks and quandle-related biquandles of these sizes. We apply the fixed point invariants of these finite biquandles to show that various virtual knots are non-trivial including the Kishino knot. We also introduce various invariants of welded knots and invariants of classical knots based on racks. Using the first mentioned invariants we produce a list of non-trivial welded knots. At least one of these has trivial quandle answering a question in \cite{FKM} and has a consequential torus with trivial group in $\R^4$, see \cite{F}. For reasons of space only the findings for $n$ up to 3 are given in an appendix. For the complete list go to www.layer8.co.uk/maths/biquandles.

We are grateful for conversations with Colin Rourke and Lou Kauffman.

\section{Notation and Definitions}
We use the notation and definitions given in \cite{FJK}, where more details can be found.

Let $X$ be a set with two binary operations, up and down, written
$$(a,b)\to a^b\hbox{ and }(a,b)\to a_b$$
The convenience of the exponential and suffix notation is that brackets can be inserted in an obvious fashion and so are not needed, For example
$$a^{bc}=(a^b)^c,\ a^{b_c}=a^{(b_c)},\ {a^b}_c=(a^b)_c, \hbox{ etc. }$$
On the other hand expressions such as $a^b_c$ are ambiguous and are not used.

Think of $b$ as acting on $a$ in both cases. We want these actions to be invertible. So there are binary operations
$$(a,b)\to a^{b^{-1}}\hbox{ and }(a,b)\to a_{b^{-1}}$$
satisfying
$$a^{bb^{-1}}=a^{b^{-1}b}=a\hbox{ and }a_{bb^{-1}}=a_{b^{-1}b}=a$$
It is convenient at this stage to introduce the function $S: X^2\to X^2$ defined by $S(a,b)=(b^a, a_b)$, known as a {\bf switch}.  

After these preliminaries we list the three axioms, B1-3 needed to define a {\bf biquandle}.

{\bf B1: } For all $a \in X$ there is a unique $x\in X$ such that $a^x=x, a=x_a$ 
and there is a unique $y\in X$ such that $a_y=y, a=y^a$ 

{\bf B2: } $S$ is invertible. So there is a function $\overline{S}: X^2\to X^2$ such that $S\overline{S}=\overline{S}S=id$

{\bf B3: } {Let $S_1=S\times id$ and $S_2=id\times S$ then $S_1S_2S_1=S_2S_1S_2$}

These axioms are a consequence of the three Reidemeister moves. If only B2 and B3 are satisfied then we have a {\bf birack}.

Axiom B1 implies that the infinite tower and the infinite well
$$x=a^{a^{a^{.^{.^{.^.}}}}},\quad y=a_{a_{a_{._{._{._.}}}}}$$
make sense.

We have $x=a_{a^{-1}}$, $y=a^{a^{-1}}$ and so B1 is equivalent to
$a^{a_{a^{-1}}}=a_{a^{-1}}\hbox{ and }a_{a^{a^{-1}}}=a^{a^{-1}}$.

In \cite{Stan} it is shown that only one half of B1 is necessary.

Using B2 write $\overline{S}(a,b)=(b_{\overline{a}},a^{\overline{b}})$.
This defines two binary operations
$$(a,b)\to a^{\overline b}\hbox{ and }(a,b)\to a_{\overline b}$$
In terms of the previous operations $\overline{a^b}$ acts down on $b$ as $a^{-1}$ and $\overline{a_b}$ acts up on $b$ as $a^{-1}$. So all these new operations are defined by the two initial up and down operations.

B3 is sometimes called the set theoretic Yang-Baxter equation. If we follow the progress of the triple $(a,b,c)$ through the two sides of the equation and swap variables we arrive at three relations true for all $a,b,c\in X$.
$$a^{c_bb^c}=a^{bc},\quad {a^b}_{c^{b_a}}={a_c}^{b_{c^a}},\quad a_{c^bb_c}=a_{bc}$$
\section{Examples}
If the down operation is trivial, so $a_b=a$ for all $a,b\in X$, then a biquandle becomes a quandle. Symmetrically, this is also the case if the up operation is trivial. A birack with trivial down (up) operation is a rack. For many examples of racks and quandles see \cite{FR}.

The simplest biquandle, with both up and down operations trivial, is the {\bf twist}. In this case $S(a,b)=(b,a)$. If $X$ has $n$ elements denote the twist by $\I_n$.

If a biquandle  has the property that the up operation (or the down operation) on its own defines a quandle then it is called a {\bf quandle-related} biquandle associated to that quandle.  In a similar fashion we can define a {\bf rack related} birack.

A biquandle is said to be {\bf linear} if it is determined by a $2\times2$ matrix $S=\pmatrix{A&B\cr C&D\cr}$, where $A, B, C, D$ are elements of an associative ring related by certain equations, see for example \cite{BuF, BuF2, FT}.

The only example of a commutative linear biquandle is given by
$$a^b=\lambda a+(1-\lambda\mu)b,\quad a_b=\mu a$$
where $a,b\in X$, a $\Z[\lambda^{\pm1}, \mu^{\pm1}]$-module \cite{Swa}. This is called the {\bf Alexander} biquandle and is denoted by $A_{\lambda\mu}(X)$. If $\mu=1$ then resulting quandle is called the {\bf Burau} quandle and is denoted by $B_\lambda(X)$. If $\mu=\lambda=1$ then we have the twist.

There are many examples of linear non-commutative biquandles.  For example let $\H$ denote the quaternion algebra with standard generators $1,i, j, k$. If $X$ is a left $\H$ module then 
$$a^b=ja+(1+i)b,\quad a_b=-ja+(1+i)b$$
is a biquandle. This is called the {\bf Budapest} biquandle and is just one of a huge family of linear non-commutative biquandles described in  \cite{F2, BF, BuF, FT}.

One example of a non-linear biquandle is due to Wada, \cite{W}. Let $X$ be a group and put 
$S(a,b)=(a^2b, b^{-1}a^{-1}b)$.

Another is due to Silver and Williams, \cite{SW}. Let $G$ be a group with a $\Z^2$ action, $a\to a\cdot x$ and $a\to a\cdot y$ where $x, y$ are a basis for $\Z^2$.  \nl
Define $S(a,b)=(b\cdot y,(b\cdot xy)^{-1}(a\cdot x)b)$.

Naturally most of the finite examples given in the appendix are non-linear and we will try to identify them when they are known.

\section{Virtual and Weld Invariants}
A biquandle is said to have {\bf order} $k$ if $k$ is the smallest positive integer for which $S^k$ is the identity. A pair of biquandles 
$(S,T)$ is a {\bf virtual} invariant if

${\bf V}$. $T$ has order two and $T_1S_2T_1=T_2S_1T_2$

A virtual invariant pair is said to be a {\bf weld} invariant if the pair satisfy

${\bf W_1}$. $T_1S_2S_1= S_2S_1T_2$

This condition is called the first {\bf forbidden} move and the pair is {\bf essential} if the second forbidden move

${\bf W_2}$. $S_1S_2T_1=T_2S_1S_2$.

is NOT satisfied; think abelian and non-abelian. Since the use of both forbidden moves unknots any virtual knot, \cite{Kan, Nel}, it is important to find essential weld invariants.

\theorem{\sl Let $S$ be a biquandle and let $T$ be the twist acting on the set $X$. Then the pair $(S,T)$ is a virtual invariant. It is a weld invariant if and only if the down operation is commutative and the down operation is trivial upstairs: in symbols $a_{bc}=a_{cb}$ and $a^{b_c}=a^b$. If $S$ is linear then it is the Alexander biquandle.}

{\bf Proof}

Condition {\bf V} is easily checked and ${\bf W_1}$ implies the two conditions.

Suppose that $S$ is linear, say $S=\pmatrix{A&B\cr C&D\cr}$. Then the second condition implies that $AD=0$. By \cite{BuF2} $D=1-A^{-1}B^{-1}AB$ so $A,B$ commute. This implies that $S$ is Alexander. \qed

There are no known linear essential weld pairs. Maybe they do not exist.

\np
\section{Finding Invariants}

Any virtual or welded knot or link can be written as the closure of a virtual braid and so can be thought of as a word $w$ in the variables $\sigma_i, \tau_i$ where $1\le i\le n-1$ and $n$ is some fixed positive integer. The variable $\sigma_i$ corresponds to a positive classical crossing and its inverse $\sigma_i^{-1}$ corresponds to a negative classical crossing. The variable $\tau_i$ corresponds to a virtual crossing or weld. 

If a pair of biquandles  $(S,T)$ is a virtual or weld invariant pair let 
$S_i=(id)^{i-1}\times S\times(id)^{n-i-1}$
and similarly for $T_i$. Suppose $w$ is a word representing some virtual braid and let $K$ be its closure.
Then the assignation
$\sigma_i\to S_i, \ \tau\to T_i$ acting on the word $w$
defines a permutation, $f_{(S,T)}$, of $X^n$  which is an invariant of the braid. 
The classic example is the {\bf Burau} map which uses the Burau quandle.

Let $\psi_1=\s_3^{-1}\s_2\s_1^2\s_2\s_4^3\s_3\s_2$, $ \psi_2=\s_4^{-1}\s_3\s_2\s_1^{-2}\s_2\s_1^2\s_2^2\s_1\s_4^5$,
$\phi_1=\s_4\s_5^{-1}\s_2^{-1}\s_1$ and $\phi_2=\s_4^{-1}\s_5^2\s_2\s_1^{-2}$.
In \cite{Big}, Bigelow defines the braids $b_1\in B_5$ and  $b_2\in B_6$ as 
$$b_1=[\psi_1^{-1}\s_4\psi_1, \psi_2^{-1}\s_4\s_3\s_2\s_1^2\s_2\s_3\s_4\psi_2]
\hbox{ and }
b_2=[\phi_1^{-1}\s_3\phi_1,\phi_2^{-1}\s_3\phi_2].$$
He goes on to show that $b_1$ and $b_2$ are non-trivial elements of the kernel of the Burau map. We will also prove that these are non-trivial by calculating  $f_{(S,T)}$ for various biquandles from the computer output.

\theorem{(Bigelow) The Burau map is not injective. In particular the braids $b_1, b_2$ defined above are non-trivial and lie in the kernel.}

{\bf Proof}\
It is possible to show that $b_1$ and $b_2$ are mapped to the identity by the Burau map, see \cite{Big}.  There are many
quandle pair representations which show that $b_1$ (having 5 strands) and $b_2$ (having 6 strands) are non-trivial.  For example, using the following virtual pair (with notation developed in the appendix)\hfil\break
$S: BQ^4_{3} \quad U = (\iota , \iota , (2 4 3) , (2 3 4)) \quad D = (\iota , (3 4) , (3 4) , (3 4))$\hfil\break
$T: BQ^4_{9} \quad U = (\iota , (3 4) , (3 4) , (3 4)) \quad D = (\iota , (3 4) , (3 4) , (3 4))$\hfil\break
the number of fixed points  of $f_{(S,T)}$ corresponding to $b_1$ is $736$. Whereas the number of fixed points of the identity permutation is of course $4^5=1024$.

The number of fixed points  of $f_{(S,T)}$ corresponding to  $b_2$ is 1648. Whereas the number of fixed points of the identity permutation is $4^6=4096$. \qed

Any $n$-tuple of $X^n$ fixed by $f_{(S,T)}$ defines a colouring of the arcs of $K$, the closure of the braid, and so the number of fixed points is also a useful invariant of $K$.
Denote the fixed points of $f_{(S,T)}$ by $F_n$. Then $f_{(S,T)}$ is invariant under the conjugation Markov move. The other Markov move corresponds to $w\leftrightarrow w\rho$ where $\rho$ is one of $\sigma_n, \sigma_n^{-1}, \tau_n$. This defines a bijection between $F_n$ and $F_{n+1}$. So a useful invariant, if $X$ is finite, is the size of $F_n$. This is referred to as the {\bf fixed-point invariant} of $K$. If $X$ is finite with $n$ elements then the  unlink with $c$ components has fixed point invariant $n^c$.

\theorem{There exists a non-trivial virtual (welded) knot which cannot be distinguished from the unknot by any linear biquandles.}

{\bf Proof}\

With the notation above let
$b=b_2\t_1\s_2\t_3\t_4\t_5\t_6$  Then the closure $K$ of $bbb$ is a knot for which the fixed-point invariant using the virtual pair
$$S: BQ^3_{3} \quad U = (\iota , (1 3 2) , (1 2 3)) \quad D = ((2 3) , (2 3) , (2 3))$$
$$T: BQ^3_{5} \quad U = ((2 3) , (2 3) , (2 3)) \quad D = ((2 3) , (2 3) , (2 3))$$
 is 9.  Here $BQ^3_{3}$ and $BQ^3_{5}$ refer to the numbering system in the appendix.
 
So $K$ is non-trivial as a virtual knot. However any linear $S$ is equivalent to the Burau representation for classical braids, see lemma 2.2 in \cite{F2}. So $b_2$ will be invisible and the sequence $\t_1\s_2\t_3\t_4\t_5\t_6$ closes to the unknot.

Using the representation determined by the essential pair
$$S: BQ^3_{3} \quad U = (\iota , (1 3 2) , (1 2 3)) \quad D = ((2 3) , (2 3) , (2 3))$$
$$T: Q^3_{1} \quad U = \iota \quad D = \iota$$
 
the fixed-point invariant for the closure of $bbb$ is 9 again, which shows that $K$ is non-trivial as a welded knot.
 \qed
 
 Finally the following application shows that the Kishino knot, $K_3$ in \cite{FJK}, is non-trivial as a virtual knot. Although this is now well known, see for example \cite{F2} and \cite{FT}, this particular knot acts as a touchstone for the effectiveness of a virtual knot invariant. The following pair
 $S: BQ^4_{53} \quad U = ((1 2 3 4) , (1 4 3 2) , (1 2 3 4) , (1 4 3 2)) \quad D = ((1 2 3 4) , (1 4 3 2) , (1 2 3 4) , (1 4 3 2))$\hfil\break
$T: BQ^4_{26} \quad U = (\iota , \iota , (1 2)(3 4) , (1 2)(3 4)) \quad D = (\iota , \iota , (1 2)(3 4) , (1 2)(3 4))$\hfil\break
results in 16 fixed points for $K_3$. Whereas the unknot would only give 4.

\section{Infinite Series from Biracks}
The {\bf writhe} of a knot diagram is the excess of positive crossings over negative crossings. It is invariant under the last two Reidemeister moves but changes by $\pm1$ under the first. Let $S$ be a birack. Suppose that a knot $K$ has a diagram with writhe $w$. Let $\phi_w$ be some scalar associated to $K$ by $S$. Define
$$\zeta_S(K)(t)=\sum^{\infty}_{w=-\infty}\phi_wt^w$$
Then $\zeta$ is an invariant of $K$ with any given writhe.

For example, let $(S,T)$ be a finite virtual pair and let  $\phi_w$  be the number of fixed points of $f_{(S,T)}$ as in the previous section. Then the invariants $\phi_w$  have the following properties.
\lemma{1. $\phi_0=|X|$ if $K$ is the unknot,

2. $\phi_1$ is the order of the set $\{x\in X|x^{x_{x^{-1}}}=x_{x^{-1}}\}$ if $K$ is the unknot,

3. $\phi_w=\phi_{-w}$.

4. If $X$ is finite then $\phi_w$ repeats for all knots. The cycle length is no greater than $|X|$ for the unknot.}

{\bf Proof}

The unknot with positive writhe $w$ can be represented as the plat closure of the braid $\sigma_1^{-w}$. In order to colour the arcs of this diagram we need to find invariant diagonal pairs of the operator $R=PQPQP\cdots$ where the number of factors is $w$ and $P, Q$ are defined by
$$P(a,b)=(b_{a^{-1}}, a^{b_{a^{-1}}}),\ Q(a,b)=(b_{a^{b^{-1}}}, a^{b^{-1}}).$$
In other words, we need the number of pairs $(a,a), a\in X$ for which $R(a,a)$ is a pair $(b,b)$ for some $b\in X$.

In \cite{FJK} the operators  $P, Q$ are denoted $S^+_-$ and $S^-_+$ respectively. They represent the transformation of the labels when the orientations of the strings are contrary and so are appropriate for plat closures.

Clearly $\phi_0=|X|$ and $\phi_1$ is the number of pairs $(x,x)$ for which both components of $P(x,x)$ are equal, that is $x^{x_{x^{-1}}}=x_{x^{-1}}$. There are similar but more complicated formul\ae\ for $\phi_w, w>1$.

If $R(x,x)=(y,y)$ then $R^{-1}(y,y)=(x,x)$ and conversely. So $\phi_w=\phi_{-w}$.

The operator $R$ defines a permutation of the subset of $X$ for which $R(x,x)$ is another diagonal pair. This permutation cannot have order greater than $|X|$.

Finally, suppose $K$ a knot or link represented as the plat closure of $\beta\in B_{2n}$. Then the appropriate operator to calculate $\phi$ is a word $R$ in $S_i^{\pm1}, P_i^{\pm1}, Q_i^{\pm1}, i\in \{1,2,\ldots,n-1\}$. Let $\Delta$ denote the 'pair' diagonal in $X^{2n}$ given by
$$\Delta=\{(x_1,x_1,x_2,x_2,\ldots,x_n,x_n)|x_i\in X\}.$$
Then $\phi_w$ for suitable $w$ is the number of elements of $\Delta$ permuted by $R$.
This permutation cannot have order greater than $|X|^n$ which is anyway finite.
Similarly $\phi_w=\phi_{-w}$ generally. 
\qed

This invariant can distinguish the trefoil from the figure eight knot as well as show they are non-trivial. It can also act as a birack invariant and distinguish some biracks of the same size. In the appendix some examples are given of this invariant.\np

\section{Notation and invariants of the output}

We will now define the notation used, look at isometries and define a few isometry invariants of biquandles.

A finite biquandle will be represented by permutations that represent the up and down actions. This will be in the form $U=(\upsilon_1,\ldots,\upsilon_n)$ and $D=(\delta_1,\ldots,\delta_n)$ where $\upsilon_i$ represents the permutation on $X_n=\{1,2,\ldots,n\}$ defined by the up operation of $i$ and $\delta_i$ represents the down operation of $i$.
The identity is written $\iota$. A whole row of identities, for example when the biquandle is a quandle, is written $I$.

As an illustration, the Alexander biquandle
$A_{12}(\Z_3)$ has
$$U=(\iota, (132), (123)),\ D=((23), (23), (23))$$
From this we can read off the operations. So if we want to find $2^1$ then we look for the up operation in $U$ corresponding to 1. This is the identity. So $2^1=2$. 
Similarly $1^2=3, 1^3=2, 2_2=3$ etc.

The twist on a set with $n$ elements. is denoted by, $\I_n$, and is defined by $\I(a,b)=(b,a)$. It is represented as a pair of identity permutations, $U=I, D=I$.

The list is defined up to isomorphism. If $\sigma$ is a permutation of $X_n$ and  $S$ is a biquandle defined by two operations, $a^b$ and $a_b$ then $\sigma$ induces an isomorphic biquandle, $\sigma_\sharp(S)$, with operations
$$a\wedge b=\sigma(\sigma^{-1}(a)^{\sigma^{-1}(b)})\hbox{ and }
a\vee b=\sigma(\sigma^{-1}(a)_{\sigma^{-1}(b)})$$
For example consider the permutation $(123)$. This changes
$$U=((23), \iota, \iota),\ D=((23), \iota, \iota)$$
into
$$U=(\iota, (13), \iota),\ D=(\iota, (13), \iota)$$

A geometric symmetry is provided by changing the sign of each crossing of a diagram. Correspondingly $S$ is replaced by $\overline S$ the inverse of $S$. Then $\overline S(a,b)=(b_{\overline a},a^{\overline b})$ defines the down operations and up operations respectively.

Another geometric symmetry is provided by changing the orientation of the knot. This swaps $U$ and $D$.

These two can be composed to get a third symmetry in which $U$ goes to $\overline D$ and $D$ goes to $\overline U$.

A quandle is a biquandle which has $D$ (or $U$) equal to the identity. Define a {\bf pseudo quandle} to be a biquandle in which all the permutations of $D$ (or $U$) are equal. A {\bf double pseudo quandle} has this property for both of $D$ and $U$. A birack is {\bf symmetric} if $U=D$.

A {\bf constant point}, $x$, of the up operation $U$ has $x^y=x^z$ for all $y,z$. Let $u$ be the number of constant points for $U$ and let $d$ be the number of constant points for $D$. Then $c_1=u+d$ and $c_2=|u-d|$ are invariant under the various symmetries. For example consider $A_{12}(\Z_3)$. Then $c_1=2,\ c_2=0$.

\section{References}

\cite{BF}  {\bf A. Bartholomew, R. Fenn,} {\it Quaternionic Invariants of Virtual Knots and Links,} JKTR
17 (2008) 231-251

\cite{Big} {\bf Stephen Bigelow} {\it The Burau representation is not faithful for n = 5} Geometry \& Topology 3 (1999) 397Ð404

\cite{BuF}  {\bf Stephen Budden, R. Fenn,} {\it Quaternion Algebras and Invariants of Virtual Knots and Links 
II: The Hyperbolic Case} JKTR 17 (2008) 279-304

\cite{BuF2} {\bf S. Budden and Roger Fenn.} {\it The equation, $[B,(A-1)(A,B)]=0$
and virtual knots and links}, Fund Math 184 (2004) pp 19-29.

\cite{F} {\bf R. Fenn,} {\it A strange torus embedding,} preprint

\cite{F2} {\bf R. Fenn,} {\it Quaternion Algebras and Invariants of Virtual Knots and Links 
I: The Elliptic Case}  JKTR 17 (2008) 279-304

\cite{FR} {\bf R. Fenn, C. Rourke,} {\it Racks and Links in Codimension Two,} JKTR, No. 4,
343-406 (1992).

\cite{FJK} {\bf R. Fenn, M. Jordan, L. Kauffman,}  {\it Biquandles and Virtual Links,} Topology and its Applications, 145 (2004) 157-175

 \cite{FKM} {\bf R. Fenn, L. Kauffman, V. Manturov,} {\it Virtual Knot Theory-Unsolved Problems }
  Fundamenta Mathematicae 188 (2005) 293-323
 
 \cite{FT} {\bf R. Fenn,  Vladimir Turaev,} {\it Weyl Algebras and Knots,}  Journal of Geometry and Physics.57 (2007) 1313-1324.

\cite{J} {\bf D.Joyce,} {\it A classifying invariant of knots; the knot quandle},
J. Pure Appl. Alg., 23 (1982) 37-65

\cite{Kan} {\bf T. Kanenobu}, {\it Forbidden moves unknot a virtual knot}, J. Knot Theory Ramifications 10 (2001), 89Ð96.

\cite{K} {\bf L.Kauffman.} {\it Virtual Knot Theory,} European J. Comb. Vol 20, 
663-690, (1999)

\cite{M}{\bf S. Matveev,} {\it Distributive groupoids in knot theory},  Math. USSR Sbornik, 47 (1984) 73-83

\cite{N} {\bf Benita Ho, Sam Nelson,} {\it Matrices and Finite Quandles,} Homology, Homotopy and Applications, Vol. 7 (2005), No. 1, pp.197-208 

\cite{Nel} {\bf S. Nelson}, {\it Unknotting virtual knots with Gauss diagram forbidden moves}, J.Knot Theory Ramifications 10 (2001), 931Ð935.

\cite{Stan} {\bf D. Stanovsky} {\it On axioms of biquandles,} J. Knot Theory Ramifications 15/7 (2006) 931--933.

\cite{SW} {\bf D.S. Silver and S.G. Williams,} {\it Alexander Groups and Virtual Links,} to
appear in JKTR.

\cite{W} {\bf M. Wada,} {\it Twisted Alexander Polynomial for finitely presentable group.} 
Topology 33 No. 2. 241-256 (1994)

 \section{Appendix: the computer output}
 In this section we will give the computer output up to size 3. For the full list go to the website
 www.layer8.co.uk/maths/biquandles

A quandle of size $n$, is indicated by $Q^n_k$ where $k$ is its position on the list of quandles of size $n$. Similarly $BQ^n_k$ denotes a biquandle which is not a quandle, $R^n_k$ a rack which is not a quandle and $BR^n_k$ a birack which is neither a biquandle or a rack. The order of the examples is determined by the computer search for biracks, which enumerates permutations of $1, \ldots, n^2$ 
to determine switches that are candidates to be biracks (axiom $B2$).  For each candidate satisfying axiom $B3$ the current list is checked for isometries before adding the new birack to the end of list.  During this process if the new birack is a rack whose representation has trivial up or down action it is always added to the list, swapping $U$ and $D$ if necessary so that the isometry with a trivial down action is added.  Any existing isometry of such a birack is removed from the list, thereby ensuring that the list always contains representations of racks with trivial down actions.  The order of the other lists is determined by the order of the list of biracks

In the list, a pseudo quandle is indicated by $PQ$, a double pseudo quandle by $DPQ$ and a symmetric one by $S$. The first example is always the twist.

These are the results for $n=2$. There is one quandle, one rack and one biquandle.

$\I_2=Q^2_{1}\quad U=I\quad D=I\quad \hbox{order 2} \quad S, c_1 = 4, c_2 = 0$\hfil\break
$BQ^2_{1}\quad U=((1 2) , (1 2))\quad D=((1 2) , (1 2))\quad \hbox{order 2} \quad S, DPQ, c_1 = 4, c_2 = 0$\hfil\break
$R^2_{1}\quad U=((1 2) , (1 2))\quad D=I\quad \hbox{order 4} \quad c_1 = 4, c_2 = 0$\hfil\break

These are the results for  $n=3$.  There are 3 quandles, 3 racks, 7 biquandles and 3 biracks, $BQ_3^3$ is $A_{12}(\Z_3)$, $BQ_5^3$ is $A_{22}(\Z_3)$ and $Q_2^3$ is $B_2(\Z_3)$.

$\I_3=Q^3_{1}\quad U=\iota\quad D=\iota\quad \hbox{order 2} \quad S, c_1 = 6, c_2 = 0$\hfil\break
$Q^3_{2}\quad U=((2 3) , \iota , \iota)\quad D=\iota\quad \hbox{order 4} \quad c_1 = 4, c_2 = 2$\hfil\break
$Q^3_{3}\quad U=((2 3) , (1 3) , (1 2))\quad D=\iota\quad \hbox{order 3} \quad c_1 = 3, c_2 = 3$\hfil\break

$R^3_{1}\quad U=((1 3 2) , (1 3 2) , (1 3 2))\quad D=\iota\quad \hbox{order 6} \quad c_1 = 6, c_2 = 0$\hfil\break
$R^3_{2}\quad U=((1 3) , \iota , (1 3))\quad D=\iota\quad \hbox{order 4} \quad c_1 = 4, c_2 = 2$\hfil\break
$R^3_{3}\quad U=((1 3) , (1 3) , (1 3))\quad D=\iota\quad \hbox{order 4} \quad c_1 = 6, c_2 = 0$\hfil\break

$BQ^3_{1}\quad U=(\iota , (2 3) , (2 3))\quad D=(\iota , (2 3) , (2 3))\quad \hbox{order 2} \quad S, c_1 = 2, c_2 = 0$\hfil\break
$BQ^3_{2}\quad U=(\iota , (2 3) , (2 3))\quad D=((2 3) , (2 3) , (2 3))\quad \hbox{order 4} \quad PQ, c_1 = 4, c_2 = 2$\hfil\break
$BQ^3_{3}\quad U=(\iota , (1 3 2) , (1 2 3))\quad D=((2 3) , (2 3) , (2 3))\quad \hbox{order 3} \quad PQ, c_1 = 3, c_2 = 3$\hfil\break
$BQ^3_{4}\quad U=(\iota , \iota , (1 2))\quad D=(\iota , \iota , (1 2))\quad \hbox{order 2} \quad S, c_1 = 2, c_2 = 0$\hfil\break
$BQ^3_{5}\quad U=((2 3) , (2 3) , (2 3))\quad D=((2 3) , (2 3) , (2 3))\quad \hbox{order 2} \quad S, DPQ, c_1 = 6, c_2 = 0$\hfil\break
$BQ^3_{6}\quad U=((1 2) , (2 3) , (1 3))\quad D=((1 2 3) , (1 2 3) , (1 2 3))\quad \hbox{order 3} \quad PQ, c_1 = 3, c_2 = 3$\hfil\break
$BQ^3_{7}\quad U=((1 2 3) , (1 2 3) , (1 2 3))\quad D=((1 3 2) , (1 3 2) , (1 3 2))\quad \hbox{order 2} \quad DPQ, c_1 = 6, c_2 = 0$\hfil\break

$BR^3_{1}\quad U=(\iota , (2 3) , (2 3))\quad D=((2 3) , \iota , \iota)\quad \hbox{order 4} \quad c_1 = 2, c_2 = 0$\hfil\break
$BR^3_{2}\quad U=((2 3) , \iota , \iota)\quad D=((2 3) , (2 3) , (2 3))\quad \hbox{order 4} \quad c_1 = 4, c_2 = 2$\hfil\break
$BR^3_{3}\quad U=((1 2 3) , (1 2 3) , (1 2 3))\quad D=((1 2 3) , (1 2 3) , (1 2 3))\quad \hbox{order 6} \quad S, c_1 = 6, c_2 = 0$\hfil\break

For $n=4$ There are  7 quandles, 12 racks, 57 biquandles and 71 biracks, the lists are available from www.layer8.co.uk/maths/biquandles.

For the cases $n=5,6$ our current programming techniques are not able to produce a full list of distinct biracks (i.e with isometries 
removed) in a reasonable time.  However we are able to produce a full list of distinct quandles.  From these lists we are able to 
produce a list of quandle-related biracks and consequently quandle-related biquandles by setting $D$ to 
be a quandle and searching for matrices $U$ for which $S=U,D$ is a birack.  In particular, since the trivial operation $I$ is a 
quandle, the search for quandle-related biracks discovers the full list of racks.  The results for the cases $n=5,6$ are summarized in 
the following table, the full lists are presented at the above website.

$$\vbox{ \offinterlineskip \halign{\strut
\vrule \hfil \quad $#$ \quad \hfil \vrule & \hfil\ #\ \hfil \vrule & \hfil\ #\ \hfil \vrule & \hfil\ #\ \hfil \vrule & \hfil\ #\ \hfil \vrule \cr 
\noalign{\hrule} 
n & quandles & racks & quandle-related biquandles & quandle-related biracks \cr 
\noalign {\hrule} 
5 & 21 & 52 & 113 & 517 \cr 
\noalign {\hrule} 
6 & 72 & 280 & 1506 & 11704 \cr 
\noalign{\hrule}}}$$

\subsection{Welded knots}

Having calculated lists of biquandles it is trivial to determine essential pairs of biquandles.  There are two such pairs of size 3:

\parindent=2em
\item{P1.} $S: BQ^3_{3} \quad U = (\iota , (1 3 2) , (1 2 3)) \quad D = ((2 3) , (2 3) , (2 3))$\hfil\break
$T: Q^3_{1}=\I_3 \quad U = \iota \quad D = \iota$\hfil\break

\item{P2.}$S: Q^3_{3} \quad U = ((2 3) , (1 3) , (1 2)) \quad D = \iota$\hfil\break
$T: BQ^3_{5} \quad U = ((2 3) , (2 3) , (2 3)) \quad D = ((2 3) , (2 3) , (2 3))$\hfil\break
\parindent=0pt

and 8 of size 4:

\parindent=2em
\item{P3.}$S: BQ^4_{3} \quad U = (\iota , \iota , (2 4 3) , (2 3 4)) \quad D = (\iota , (3 4) , (3 4) , (3 4))$\hfil\break
$T: Q^4_{1} =\I_4\quad U = \iota \quad D = \iota$\hfil\break

\item{P4.}$S: BQ^4_{19} \quad U = (\iota , (1 3)(2 4) , (1 4)(2 3) , (1 2)(3 4)) \quad D = ((2 4 3) , (2 4 3) , (2 4 3) , (2 4 3))$\hfil\break
$T: Q^4_{1} =\I_4\quad U = \iota \quad D = \iota$\hfil\break

\item{P5.}$S: BQ^4_{34} \quad U = (\iota , \iota , \iota , (1 3 2)) \quad D = ((2 3) , (1 3) , (1 2) , (1 2 3))$\hfil\break
$T: BQ^4_{23} \quad U = (\iota , \iota , \iota , (1 3 2)) \quad D = (\iota , \iota , \iota , (1 2 3))$\hfil\break

\item{P6.}$S: BQ^4_{38} \quad U = ((2 3 4) , (2 3 4) , (2 3 4) , (2 3 4)) \quad D = ((2 4 3) , (2 3) , (3 4) , (2 4))$\hfil\break
$T: BQ^4_{39} \quad U = ((2 3 4) , (2 3 4) , (2 3 4) , (2 3 4)) \quad D = ((2 4 3) , (2 4 3) , (2 4 3) , (2 4 3))$\hfil\break

\item{P7.}$S: BQ^4_{41} \quad U = ((2 3 4) , (1 3 2) , (1 4 3) , (1 2 4)) \quad D = ((2 3 4) , (2 3 4) , (2 3 4) , (2 3 4))$\hfil\break
$T: BQ^4_{39} \quad U = ((2 3 4) , (2 3 4) , (2 3 4) , (2 3 4)) \quad D = ((2 4 3) , (2 4 3) , (2 4 3) , (2 4 3))$\hfil\break

\item{P8.}$S: BQ^4_{56} \quad U = ((1 2)(3 4) , (1 2)(3 4) , (1 2)(3 4) , (1 2)(3 4)) \quad D = ((1 3 2) , (1 2 4) , (1 4 3) , (2 3 4))$\hfil\break
$T: BQ^4_{50} \quad U = ((1 2)(3 4) , (1 2)(3 4) , (1 2)(3 4) , (1 2)(3 4)) \quad D = ((1 2)(3 4) , (1 2)(3 4) , (1 2)(3 4) , (1 2)(3 4))$\hfil\break

\item{P9.}$S: Q^4_{6} \quad U = ((2 4) , (1 3) , (2 4) , (1 3)) \quad D = \iota$\hfil\break
$T: BQ^4_{50} \quad U = ((1 2)(3 4) , (1 2)(3 4) , (1 2)(3 4) , (1 2)(3 4)) \quad D = ((1 2)(3 4) , (1 2)(3 4) , (1 2)(3 4) , (1 2)(3 4))$\hfil\break

\item{P10.}$S: BQ^4_{51} \quad U = ((1 2)(3 4) , (1 3)(2 4) , (1 3)(2 4) , (1 2)(3 4)) \quad D = ((1 2 4 3) , (1 2 4 3) , (1 2 4 3) , (1 2 4 3))$\hfil\break
$T: Q^4_{1} =\I_4\quad U = \iota \quad D = \iota$\hfil\break
\parindent=0pt

The website also lists the 17 essential pairs of quandle-related biquandles of size 5 and the 271 essential pairs of quandle-related 
biquandles of size 6.  However, we include three essential pairs of size 6 here for reference below.

\parindent=2em
\item{P11.}$S: BQ^6_{10} \quad U = (\iota , (1 3 4 5 6) , \iota , \iota , \iota , \iota) \quad D = ((3 4 6 5) , (1 6 5 4 3) , (1 6 4 5) , (1 5 6 3) , (1 4 3 6) , (1 3 5 4))$\hfil\break
$T: BQ^6_{1494} \quad U = (\iota , (1 3 4 5 6) , \iota , \iota , \iota , \iota) \quad D = (\iota , (1 6 5 4 3) , \iota , \iota , \iota , \iota)$\hfil\break

\item{P12.}$S: BQ^6_{22} \quad U = (\iota , (1 3 4 5 6) , \iota , \iota , \iota , \iota) \quad D = ((3 6)(4 5) , (1 6 5 4 3) , (1 4)(5 6) , (1 6)(3 5) , (1 3)(4 6) , (1 5)(3 4))$\hfil\break
$T: BQ^6_{1494} \quad U = (\iota , (1 3 4 5 6) , \iota , \iota , \iota , \iota) \quad D = (\iota , (1 6 5 4 3) , \iota , \iota , \iota , \iota)$\hfil\break

\item{P13.}$S: BQ^6_{49} \quad U = ((3 4 5 6) , \iota , \iota , \iota , \iota , \iota) \quad D = ((3 6 5 4) , (3 4 5 6) , (2 4 6 5) , (2 5 3 6) , (2 6 4 3) , (2 3 5 4))$\hfil\break
$T: BQ^6_{230} \quad U = ((3 4 5 6) , \iota , \iota , \iota , \iota , i) \quad D = ((3 6 5 4) , i , i , i , i , i)$\hfil\break
\parindent=0pt

Given that the fixed-point invariant $F_n$ takes the value $n$ for the unknot, it is possible to use the list of essential pairs to 
search for distinct non-trivial welded knots.  Using this technique the following knots have been distinguished.

$$\vbox{ \offinterlineskip \halign{\strut
\vrule \hfil \quad # \quad \hfil \vrule & \hfil\ #\ \hfil \vrule & \hfil\ #\ \hfil \vrule  & \hfil\ #\ \hfil \vrule  
& \hfil\ #\ \hfil \vrule & \hfil\ #\ \hfil \vrule & \hfil\ #\ \hfil \vrule & \hfil\ #\ \hfil \vrule & \hfil\ #\ \hfil \vrule \cr 
\noalign{\hrule} 
& & \multispan5 \hfil Size of $F_n$ for pair \hfil \vrule \cr
\omit\vrule & \omit\vrule & \omit\hrulefill\vrule & \omit\hrulefill\vrule & \omit\hrulefill\vrule & \omit\hrulefill\vrule & \omit\hrulefill \cr
\omit\vrule height2pt \hfil \vrule height2pt&                 &     &     &     &     &      \cr
 knot & braid word                      & P3  & P4  & P11 & P12 & P13  \cr 
\noalign {\hrule} 
 w3.1 & $\s1\t2\s3-\s2-\s2-\s1\t2-\s3\s2 $             & 10 & 4  & 6  & 6  & 6   \cr 
 w3.2 & $\t1-\s2\t1-\s1-\s1\t2           $             & 10 & 16 & 6  & 6  & 6   \cr
 w4.1 & $\s1\t1-\s1\s2\s1\t1-\s1-\s2      $            & 10 & 4  & 26 & 6  & 6   \cr
 w4.2 & $-\s1-\s2\s3\t2\s1-\s4\s3\t2\s3\s4-\s3-\s2 $   & 10 & 4  & 6  & 6  & 26  \cr
 w4.3 & $-\s1\s2\s3\t2\s1-\s4\s3\t2\s3\s4-\s3\s2    $  & 4  & 4  & 26 & 26 & 6   \cr
 w4.4 & $-\s1\s2\s3\t2\s1-\s4\s3-\s2\s3\s4-\s3\t2  $   & 4  & 4  & 6  & 26 & 6   \cr
 w4.5 & $\t1\s2-\s1\t1\s1\s2                  $        & 4  & 16 & 6  & 6  & 6   \cr
 w4.6 & $-\s1-\s2\t3-\s2\s1-\s4\t3-\s2-\s3\s4-\s3\s2 $ & 4  & 4  & 6  & 26 & 26  \cr
 w6.1 & $-\s1-\s2-\s2-\s2\s1-\s3-\s2-\s2-\s2\s3\t2 $   & 28 & 28 &    &    &     \cr
\noalign{\hrule}}}$$
The non-trivial welded knot w4.5 has trivial fundamental quandle and so trivial fundamental group, \cite{F}. This answers a question in \cite{FKM}.

\subsection{Birack Series}

The series of coefficients described in section 6 can be used to distinguish classical and virtual knots from the unknot. We append here a table of coefficients which generate the repeating invariants for the unknot, the trefoil, the figure eight and two of the Kishino knots, see \cite{FJK}.  The virtual pairs $(S,T)$ were constructed with $S$ chosen to be a rack or birack, shown in the table using the above numbering scheme, and $T$ set to the twist in each case.  Many such examples exist, note in fact that the column corresponding to $R^5_{40}$ is redundant.

$$\vbox{ \offinterlineskip \halign{\strut
\vrule \hfil \quad # \quad \hfil \vrule & \hfil\ #\ \hfil \vrule & \hfil\ #\ \hfil \vrule  & \hfil\ #\ \hfil \vrule  
& \hfil\ #\ \hfil \vrule & \hfil\ #\ \hfil \vrule  & \hfil\ #\ \hfil \vrule \cr 
\noalign{\hrule} 
& & \multispan3 \hfil Coefficient series for virtual pair constructed from\hfil \vrule \cr
\omit\vrule & \omit\vrule & \omit\hrulefill\vrule & \omit\hrulefill\vrule & \omit\hrulefill \cr
\omit\vrule height2pt \hfil \vrule height2pt&                 &     &     & \cr
 knot & braid word         & $R^5_{40}$  & $R^6_{114}$ & $BR^6_{125}$ \cr 
\noalign {\hrule} 
     & unknot              &  5 3      &  6 4        &  6 3 3     \cr
 3.1 & $\s_1\s_1\s_1 $             &  11 9   &  18 16  &  12 9 9   \cr 
 4.1 & $\s_1\s^{-1}_2\s_1\s^{-1}_2 $         &  5 3     &  18 16  &  6 3 3     \cr
 K1  & $\s_1\s^{-1}_2\s^{-1}_1\t_2\s_1\s_2\s^{-1}_1\t_2$ &  11 9  &  6 4       &  6 3 3     \cr
 K2  & $\s^{-1}_1\s^{-1}_2\s_1\t_2\s^{-1}_1\s_2\s_1\t_2$ &  5 3      &  6 4         &  12 9 9   \cr
\noalign{\hrule}}}$$

\bye